 \documentclass[12pt]{article}
 \usepackage{amsmath,amsfonts,theorem,epsfig}
 \numberwithin{equation}{section}
 \newtheorem{prop}{Proposition}[section]
 \newtheorem{con}[prop]{Conjecture}
 \newtheorem{cor}{Corollary}[section]
 \theorembodyfont{\rmfamily}
 \newtheorem{dfn}[prop]{Definition}
 \newcommand{\half}{\frac12}
 \DeclareMathOperator{\Ad}{Ad}
 \DeclareMathOperator{\id}{id}
 \DeclareMathOperator{\Tr}{Tr}
 \DeclareMathOperator{\CL}{CL}
 \DeclareMathOperator{\Rep}{Rep}
 \DeclareMathOperator{\CLRep}{CLRep}
 \newcommand{\conj}{\mathrm{conj}}
  \DeclareMathOperator{\Conj}{Conj}

 \newcommand{\rest}[1]{{}_{{\textstyle|}#1}}
 \newcommand{\p}{\partial}
 \newcommand{\wave}{\widetilde}
 \newcommand{\we}{\wave{e}}
 \newcommand{\wg}{\wave{g}}
 \newcommand{\wu}{\wave{u}}
 \newcommand{\wv}{\wave{v}}
 \newcommand{\wGa}{\wave{\Ga}}
 \newcommand{\wTh}{\wave{\Theta}}
 \newcommand{\sA}{\mathcal A}
 \newcommand{\sB}{\mathcal B}
 \newcommand{\sG}{\mathcal G} % sheaf G
 \newcommand{\al}{\alpha}
 \newcommand{\De}{\Delta}
 \newcommand{\Ga}{\Gamma}
 \newcommand{\Si}{\Sigma}
 \newcommand{\Z}{\mathbb Z}

 \DeclareMathOperator{\res}{res}
 \DeclareMathOperator{\Hom}{Hom}

 \DeclareMathOperator{\im}{im}

 \newcommand{\U}{\operatorname U} % unitary group
 \newcommand{\SU}{\operatorname{SU}}

\begin{document}

 \title{Lattice gauge theories and \\ the Florentino conjecture}
 \markright{\hfill Lattice gauge theories and the Florentino conjecture
\quad}

 \author{Andrei Tyurin}
 \date{Feb 2001 }
 \maketitle

 \begin{abstract}
We study the relation between the space of representation classes of the
fundamental group of a Riemann surface and gauge theory on trivalent
graphs. We construct a partial gauge fixing in the latter gauge theory. As
an application we get a proof of a conjecture of Florentino.
 \end{abstract}

\section{Introduction}
\label{s1}

Two boundary cases of gauge theories are extremely useful for
explorations: gauge theory with a discrete gauge group (see Freed and
Quinn \cite{FQ}) and gauge theory over a discrete space, usually called a
``lattice model''. Pragmatically, such lattice models are required as a
test of the concepts. In this paper we explore from a ``theoretical''
point of view the second case: gauge theory on graphs, more precisely, on
trivalent graphs.

At the end of \cite[Seminar~6]{1}, Michael Atiyah proposed a programme to
make sense of the functional integral of Chern--Simons theory for finite
level $k$ in a rigorous way. Because direct attempt at the analysis is
extremely difficult he proposed to think of connections in purely
combinatorial terms. In particular using a discrete analogue of gauge
theory, one should express the Chern--Simons functional in a combinatorial
framework. For this, we investigated geometric properties of the phase map
in \cite{2}. On the other hand, we extended to spin networks all the
constructions of TQFT concerning Wilson loops (knots). This made possible
a partial realization of Atiyah's programme using trivalent graphs instead
of loops.

Unfortunately we don't have enough space to recall all the details of this
treatment. See \cite{3} and \cite{4}, from which the current paper derives.

Our starting point is of course exactly the same as in a number of texts
about CFT, quantum gravity and spin networks. Our previous paper \cite{5}
contains a brief survey of mathematical aspects of $\SU(2)$-spin networks
of genus $g$. Here we extend the study of the ``pumping trick''
geometry to obtain a close relation between gauge theory on graphs and
gauge theory on Riemann surfaces. The construction and results can be
summarized as follow. Any graph $\Ga$ defines a handlebody $\wGa$, that
is, a 3-manifold with boundary a Riemann surface $\p\wGa=\Si_{\Ga}$ by
pumping up the edges of $\Ga$ to tubes and the vertices to small
2-spheres with 3 holes. To investigate the relation between the space of
classes of $\SU(2)$-representations of the fundamental group of the
Riemann surface $\p\wGa=\Si_{\Ga}$ and the $\SU(2)$-gauge theory on
$\Ga$ we construct the partial gauge fixing in the latter gauge theory.
As an application of this gauge fixing we get the proof of the Florentino
conjecture.

To be concrete we only consider trivalent graphs and start by recalling
some foundations of graph geometry in a form that is convenient for us.

\section{Geometry of trivalent graphs}

 Write $E(\Ga)$ and $V(\Ga)$ for the set of edges and vertices of a graph
$\Ga$, and $E(\Ga)_v$ for the set of edges out of $v$. Write $F(\Ga)=\{v\in
e\}$ and $L(\Ga)$ for the set of flags and loops of $\Ga$, where a flag is
an edge with a fixed end, and a loop is an edge with only one vertex. The
latter set is a subset $L(\Ga)\subset E(\Ga)$ and sending a loop to its
vertex gives the map $v\colon L(\Ga)\to V(\Ga)$ which is an embedding:
\begin{equation}
L(\Ga)=v(L(\Ga)) \subset V(\Ga)
\end{equation}
because we only consider trivalent graphs.

The two projections $ e\colon F(\Ga)\to E(\Ga)$ and $v\colon F(\Ga)\to
V(\Ga)$ are ramified covers of degree 2 and 3 having the same ramification
locus $W_e=W_v \subset F(\Ga)$, consisting of pairs $v\in e$ where the
edge $e$ is a loop. The branch loci of these covers are
\begin{equation}
R_e=L(\Ga) \subset E(\Ga) \quad\hbox{and} \quad R_v=v(L(\Ga)) \subset
V(\Ga).
\end{equation}
Write $\vert\ \vert$ for the number of elements of a finite set. Then
 \begin{equation}
2 \cdot \vert E(\Ga) \vert-\vert L(\Ga) \vert=\vert F(\Ga) \vert=
3 \cdot \vert V(\Ga) \vert-\vert L(\Ga)\vert.
 \end{equation}

Hence $\vert V(\Ga) \vert=2g-2$ and $ \vert E(\Ga) \vert=3g-3$,
where $g>1$ is a certain integer called the genus of $\Ga$.

A path of length 1 in $\Ga$ is just an oriented edge $\vec e$. Let
$P_1(\Ga)$ be the set of paths of length 1 in $\Ga$. Reversing the
orientation of an edge $e$ defines an involution
\begin{equation}
i_e\colon P_1(\Ga)\to P_1(\Ga).
 \end{equation}
These edgewise involutions generate a group $O(\Ga) =\Z_2^{3g-3}$ acting
on $P_1(\Ga)$, with quotient
\begin{equation}
P_1(\Ga) / O(\Ga)=E(\Ga),
\end{equation}
and the corresponding map
\begin{equation}
O\colon P_1(\Ga)\to E(\Ga).
\end{equation}
Any section
\begin{equation}
o\colon E(\Ga)\to P_1(\Ga)
\end{equation}
of this map is called an orientation of $\Ga$. The group $O(\Ga)$ contains
a diagonal involution
\begin{equation}
i_\De\colon P_1(\Ga)\to P_1(\Ga)
 \end{equation}
reversing the orientation of all 1-paths.

Around every vertex $v\in V(\Ga)$ there are three local edges
$e_1,e_2,e_3$ or three local paths $\vec e(v)_1,\vec e(v)_2,\vec e(v)_3$
(local means that a priori two of them can be part of a loop) or local
flags. Let $P_1^v$ be the set of 1-paths with a vertex $v\in V(\Ga)$. The
orientation forgetting projection $f$ sends $P_1^v$ to $E(\Ga)_v$.

Every path $\vec e\in P_1(\Ga)$ defines two vertices
\begin{equation}
v_s(\vec e), \quad v_t(\vec e)\in V(\Ga)
\end{equation}
the target and the source ends of an arrow.
Of course
 \[
e\in L(\Ga) \implies v_s(\vec e)=v_t(\vec e) .
 \]
but even such vertex admits two local flags $\vec e_s \neq \vec e_t$.

\begin{dfn} A section $o$ is called a {\em decay orientation} if every
vertex $v\in V(\Ga)$ is a target and a source for its local edges.
\end{dfn}
More precisely, any orientation $o$ defines a section of the projection
$f\colon P_1^v\to E(\Ga)_v$ for any vertex $v$. This section must have two
arrows with opposite orientation.

We have maps
 \begin{equation}
 \begin{matrix}
v_s\colon P_1(\Ga)\to V(\Ga), \quad
 v_t\colon P_1(\Ga)\to V(\Ga) \\[8pt]
(v_s \times v_t)\colon P_1(\Ga)\to V(\Ga) \times V(\Ga)
 \end{matrix}
 \end{equation}
and the intersection of the image with the diagonal
\begin{equation}
(v_s \times v_t) ( P_1(\Ga)) \cap V(\Ga)_\Delta=v(L(\Ga)).
\end{equation}

A path of length $d$ in $\Ga$ is an ordered sequence $\vec e_1,\dots,
\vec e_d$ of oriented edges such that for every number $i$
\begin{equation}
v_t(\vec e_i)=v_s(\vec e_{i+1}).
\end{equation}
Let $P_d(\Ga)$ be the set of paths of length $d$ in $\Ga$. Every path
$(\vec e_1,\dots, \vec e_d)\in P_d(\Ga)$ defines two vertices
\begin{equation}
v_s(\vec e_1,\dots, \vec e_d), \quad\hbox{and}\quad v_t(\vec e_1,\dots,
\vec e_d)\in V(\Ga)
\end{equation}
the source and target of a path; this defines maps
 \begin{equation}
 \begin{matrix}
v_s\colon P_d(\Ga)\to V(\Ga), \quad
v_t\colon P_d(\Ga)\to V(\Ga) \\[8pt]
(v_s \times v_t)\colon P_d(\Ga)\to V(\Ga) \times V(\Ga).
 \end{matrix}
 \end{equation}
The inverse image of the intersection with the diagonal
\begin{equation}
L_d(\Ga)=(v_s \times v_t)^{-1} ( V(\Ga)_\Delta )
\end{equation}
is the set of loops of length $d$. In particular,
 \[
L_1(\Ga)= L(\Ga).
 \]
For a vertex $v\in V(\Ga)$ we have the set of loops
\begin{equation}
L_d(\Ga)_v=(v_s \times v_t)^{-1} ( v)
\end{equation}
and the union
\begin{equation}
L_\infty(\Ga)_v=\bigcup_{d=1}^\infty L_d(\Ga)_v
\end{equation}
admits a group structure. This is the combinatorial fundamental group
\begin{equation}
\pi_1^C(\Ga)_v=L_\infty(\Ga)_v
\end{equation}
of our graph $\Ga$. It is easy to check

\begin{prop} \begin{enumerate} \item The combinatorial fundamental group
$\pi_1^C(\Ga, v_0)$ is equal to the usual fundamental group $\pi_1(\Ga,
v_0)$,
\item the usual fundamental group $\pi_1(\Ga)$ is the free group with
$g$ generators, where $g$ is genus of $\Ga$.
\end{enumerate}
\end{prop}
Indeed, contracting $2g-3$ simple edges of $\Ga$ gives a bouquet of $g$
circles.

Now consider the $\Z$-module $V^{\Z}$ of all formal linear combinations of
vertices of $V(\Ga)$ with coefficients in $\Z$. Then any graph $\Ga$ of
genus $g$ can be represented as an endomorphism (matrix)
\begin{equation}
q_{\Ga}\colon V^{\Z}\to V^{\Z}
\end{equation}
with coefficient $\al_{v_i,v_j}=$ the number of edges joining $v_i$ and
$v_j$. Thus every graph $\Ga$ defines an integer quadratic form $q_{\Ga}$
on the $\Z$-module $V^{\Z}$. For example,
\begin{align}
\hbox{for the genus 2 loop free graph $\Theta$:}\quad
& q_\Theta=\begin{pmatrix} 0&3\\ 3&0 \end{pmatrix} \\[8pt]
\hbox{for the 2-loop graph of genus 2:} \quad
& q_{\Ga_{2}}=\begin{pmatrix} 2&1\\ 1&2 \end{pmatrix} \\[8pt]
\hbox{for the genus 3 multi-theta graph:} \quad
& q_{\Theta_{3}}=\begin{pmatrix}
0&2&1&0\\
2&0&0&1 \\
1&0&0&2 \\
0&1&2&0
\end{pmatrix}
\end{align}
(see Figure~1, p.~\pageref{fig!1}). All the coefficients are non\-negative,
and for a trivalent graph, the sum of the coefficients along every row and
column is equal to 3.

The permutation group $S_{2g-2}$ acts on $V_{2g-2}$ by renumbering the
vertices, and transforms matrices of the quadratic forms.

\begin{dfn} As usual, a graph $\Ga$ is {\em hyperbolic} if
there are two subsets $V_+,V_-\subset V_{2g-2}$ such that the subspaces
$V_{\pm}^{\Z}$ are isotropic with respect to $q_{\Ga}$.
\end{dfn}

For such a graph the matrix $q_{\Ga}$ has the block form
 \begin{equation}
 q_{\Ga}=\begin{pmatrix}
0&0&*&*\\
0&0&*&*\\
*&*&0&0 \\
*&*&0&0
\end{pmatrix}
 \end{equation}
where the blocks
\begin{equation}
\half q_{\Ga}=\begin{pmatrix}
 *&*\\
 *&*\\
\end{pmatrix} \in \Hom(V_+^{\Z}, V_-^{\Z})
 \end{equation}
satisfy the same conditions on rows and columns as before. In this case
the set of edges $E(\Ga)$ can be presented as a collection of triples
 \begin{equation}
E(\Ga)=\bigcup_{v\in V_+} E(\Ga)_v
 \end{equation}
and we have the map
\begin{equation}
v_+\colon E(\Ga)\to V_+
\end{equation}
which defines an orientation
\begin{equation}
o_{V_+}\colon E(\Ga)\to P_1(\Ga).
\end{equation}
If our graph $\Ga$ is connected then this is a decay orientation
(see Definition~\ref{def2.1}).

Recall that the $\Z$-module of formal linear combimations over $\Z$ of
trivalent graphs is an algebra with respect to the multiplication induced
by disjoint union of graphs. For example,
\begin{equation}
 q_{\Theta^{2}}=\begin{pmatrix}
0&3&0&0\\
3&0&0&0 \\
0&0&0&3 \\
0&0&3&0
\end{pmatrix}
 \end{equation}
and so on.

\section{$\SU(2)$-gauge theory on trivalent graphs}

\subsection{Space of connections}

A connection $a$ on the trivial $\SU(2)$-bundle on $\Ga$ is a map
\begin{equation}
a\colon P_1(\Ga)\to \SU(2)
\end{equation}
such that for the involution (2.5) we have
\begin{equation}
a(i_e(\vec e))=a(\vec e)^{-1}.
\end{equation}
Then the ``path integral'' is given by
\begin{equation}
a(\vec e_1,\dots, \vec e_d)=a(\vec e_1) \cdots a(\vec e_d)\in
\SU(2),
\end{equation}
that is, it defines a map
\begin{equation}
a\colon P_d(\Ga)\to \SU(2)
\end{equation}
such that for the orientation reversing involution $i_\De$ we have
\begin{equation}
a(i(\vec e_1,\dots, \vec e_d))=a(\vec e_1,\dots, \vec e_d)^{-1}.
\end{equation}
In the same vein we have the monodromy map for loops
\begin{equation}
a\colon L_d(\Ga)\to \SU(2).
\end{equation}
Obviously every connection is {\it flat\/}.

Let $\sA(\Ga)$ be the space of connections, that is, the space of
functions (3.1) subject to the constraint (3.2):
\begin{equation}
\sA(\Ga)=\bigl\{ a\in \SU(2)^{P_1(\Ga)} \bigm| a(i_e(\vec e))=a(\vec
e)^{-1} \bigr\}.
\end{equation}

\subsection{Gauge transformations group}

Every element $\wg$ of the gauge transformations group
$\sG(\Ga)$ is a function
\begin{equation}
\wg\colon V(\Ga)\to \SU(2),
\end{equation}
that is,
\begin{equation}
\sG(\Ga)=\SU(2)^{V(\Ga)}
\end{equation}
with componentwise multiplication. This group acts on the space of
connections $\sA(\Ga)$ by the rule
\begin{equation}
\wg(a(\vec e))=\wg(v_s(\vec e)) \cdot a(\vec e)
\cdot
\wg(v_t(i(\vec e))).
\end{equation}
Recall that $\wg(v_t(i(\vec e)))= \wg(v_t(\vec
e))^{-1}.$

The space of gauge orbits
\begin{equation}
\sB(\Ga)=\sA(\Ga) /\sG(\Ga)
\end{equation}
is the space of representation classes of the fundamental group of $\Ga$,
as expected, because all our connections are flat:
\begin{equation}
 \sA(\Ga) /\sG(\Ga)=\CLRep(\pi_1(\Ga)),
\end{equation}
where as usual for a group $G$
\begin{equation}
\Rep(G)=\Hom( G, \SU(2))
\end{equation}
 is the space of representations and the quotient by the adjoint action
\begin{equation}
 \CLRep(G)=\Hom( G, \SU(2)) /\Ad \SU(2)
\end{equation}
is the space of representation classes.

Now the gauge transformation group $\sG(\Ga)$ contains the diagonal
subgroup
\begin{equation}
\SU(2)_\De \subset \sG(\Ga)
\end{equation}
of constant functions (3.8). The action of this subgroup defines
the space of constant orbits
\begin{equation}
 \CL \sA(\Ga)=\sA(\Ga) / \SU(2)_\De.
\end{equation}

For every vertex $v\in V(\Ga)$ we have the subgroup
\begin{equation}
\sG(\Ga)_v=\bigl\{\wg\in \sG(\Ga) \bigm| \wg(v)=\id \bigr\}
\end{equation}
of gauge transformations preserving the framing at $v$. This is a normal
subgroup, and
\begin{equation}
\sG(\Ga) /\sG(\Ga)_v=\SU(2)_\De.
\end{equation}
Thus the full group of gauge transformations is the semidirect product
of $\sG(\Ga)_v$ and $\SU(2)_\De$.

The quotient
\begin{equation}
\sA(\Ga) /\sG(\Ga)_v=\Rep(\pi_1(\Ga, v))
\end{equation}
is the orbit space of connections framed at $v$; this space depends on
the choice of $v$.

Thus the quotient map (3.12)
\begin{equation}
P\colon \sA(\Ga)\to \CLRep(\pi_1(\Ga))
\end{equation}
can be decomposed as follows
\begin{equation}
\sA(\Ga) \xrightarrow{\,P\,} \Rep(\pi_1(\Ga)) \xrightarrow{/\SU(2)_\De}
\CLRep(\pi_1(\Ga))
\end{equation}
or
\begin{equation}
\sA(\Ga)\xrightarrow{/\SU(2)_\De} \CL\sA(\Ga) \xrightarrow{P_{cl}}
\CLRep(\pi_1(\Ga)).
\end{equation}
The involution $i_\De$ acts on the space $\sA(\Ga)$ of connections
\begin{equation}
i^*_\De\colon \sA(\Ga)\to \sA(\Ga)
\end{equation}
and there exists an element
\begin{equation}
\wg_i=\begin{pmatrix}0&1\\-1&0\end{pmatrix}
\in\SU(2)_\De \subset \sG(\Ga)
\end{equation}
such that
\begin{equation}
i^*_\De=\wg_i.
\end{equation}
(Here $i=\sqrt{-1}$.) Thus the involution $i^*_\De$ (3.23) acts trivially
on $\CL\sA(\Ga)$. Recall that a gauge fixing is a section of the
projection $P$ (3.20).

\begin{dfn}\label{def2.1} A partial gauge fixing is a section
 \[
s\colon \CLRep(\pi_1(\Ga))\to \CL\sA(\Ga).
 \]
\end{dfn}

\subsection{The map $\conj$}

Recall that the set $\Conj (\SU(2))$ of conjugacy classes of elements of
$\SU(2)$ can be described as the interval $[0,1]$ with the map
\begin{equation}
\conj\colon \SU(2)\to \Conj (\SU(2))=[0,1]
\end{equation}
that sends a matrix $g\in \SU(2)$ to
\begin{equation}
 \conj g=\frac{1}{\pi} \cdot \cos^{-1}\Bigl(\frac{1}{2}\Tr g \Bigr)\in
 [0,1].
 \end{equation}

Using this map coordinatewise gives the map
\begin{equation}
\conj\colon \CL\sA(\Ga)\to [0,1]^{3g-3}=\prod_{e\in E(\Ga)}
[0,1]_e,
\end{equation}
which is obviously surjective; it is the composite
\begin{equation}
 \sA(\Ga) \xrightarrow{\SU(2)_\De} \CL\sA(\Ga)\to
[0,1]^{3g-3}=\prod_{e\in E(\Ga)} [0,1]_e,
\end{equation}
and the involution $i_\De$ preserves its fibers.

 \subsection{Abelian gauge theory}
 Abelian gauge theory on $\Ga$ is a good model for non-Abelian gauge
theory. A $\U(1)$-connection $a$ is a map
\begin{equation}
a\colon P_1(\Ga)\to \U(1)
\end{equation}
such that
\begin{equation}
a(i_{e}(\vec e))=a(\vec e)^{-1}.
\end{equation}
The ``path integral'' is given by
\begin{equation}
a(\vec e_1,\dots, \vec e_d)=\prod_{\vec e\in P_1(\Ga)} a(\vec e)\in
\U(1),
\end{equation}
and so on. Let
\begin{equation}
\sA_{\U(1)}(\Ga)=\bigl\{ a\in \U(1)^{P_1(\Ga)} \bigm| a(i_e(\vec e))=a(\vec
e)^{-1} \bigr\}.
\end{equation}
be the space of $\U(1)$-connections. Then we have the same involution
\begin{equation}
 \sA_{\U(1)}(\Ga) \xrightarrow{i^*_e} \sA_{\U(1)}(\Ga).
\end{equation}
Every element $\wu$ of the gauge transformations group
$\sG_{\U(1)}(\Ga)$ is a function
\begin{equation}
\wu\colon V(\Ga)\to \U(1).
\end{equation}
This group acts on the space of connections $\sA_{\U(1)}(\Ga)$ by the same
rule:
\begin{equation}
\wu(a(\vec e))=\wu(v_s(\vec e)) \cdot a(\vec e)
\cdot
\wu(v_t(i(\vec e))).
\end{equation}

There is however one important difference between the Abelian and
non-Abelian theories. Namely the diagonal group
\begin{equation}
\U(1)_\De \subset \sG_{\U(1)}(\Ga)
\end{equation}
acts trivially. Thus
\begin{equation}
\sB_{\U(1)}(\Ga)=\sA_{\U(1)}(\Ga) /\sG_{\U(1)}(\Ga)=\Hom(\pi_1(\Ga) ,
\U(1))=\U(1)^g.
\end{equation}
Let
\begin{equation}
P_a\colon \sA_{\U(1)}(\Ga)\to \U(1)^g
\end{equation}
be the projection map.

The Abelian and non-Abelian theories are related by the following map:
\begin{equation}
 d\colon \sA_{\U(1)}(\Ga)\to \sA(\Ga)
 \quad\hbox{given by}\quad
 d(a(\vec e))=\begin{pmatrix} e^{i \phi}& 0\\ 0& e^{- i \phi}
\end{pmatrix}.
\end{equation}
$d$ is equivariant with respect to every involution $i^*_e$.
Obviously the image $d (\sA_{\U(1)}(\Ga) ) \subset \sA(\Ga)$
is a 2-section of the projection $\conj$, that is, the composite
\begin{equation}
 \conj \circ d\colon \sA_{\U(1)}(\Ga)\to \prod_{e\in E(\Ga)}
[0,1]_e
\end{equation}
 is the factorization by the involution $i^*_\De$.

It is easy to check

\begin{prop} In the chain of maps
\begin{equation}
\sA_{\U(1)}(\Ga) \xrightarrow{\,d\,} \sA(\Ga) \xrightarrow{\,P\,}
\U(1)^g\in \CLRep (\pi_1(\Ga))
\end{equation}
the composite $d \circ P$ is the projection map $P_a$ (3.39).
\end{prop}

\section{From trivalent graph to handlebody}
\subsection{Pumping}
Any trivalent graph $\Ga$ defines a handlebody $\wGa$, that is, a
3-manifold with boundary $\p \wGa=\Si_{\Ga}$ by the ``pumping trick'' (see
\cite{5}): pump up the edges of $\Ga$ to tubes and the vertices to small
2-spheres. This gives a Riemann surface $\Si_{\Ga}$ of genus $g$ with a
tube $\we$ for every $e\in E(\Ga)$ and a {\it trinion} $\wv$ for every
$v\in V(\Ga)$, where each trinion is a 2-sphere with three disjoint holes.
The isotopy classes of meridian circles of tubes define $3g-3$ disjoint,
non\-contractible, pairwise nonisotopic circles $\{C_e\}, e\in E(\Ga)$
on $\Si$. The complement
 is the union
\begin{equation}
 \Si_{\Ga} \setminus \bigcup_{e\in E(\Ga)}C_e=
 \bigcup_{i=1}^{2g-2} \wv_i
 \end{equation}
 of $2g-2$ trinions (pairs of pants) corresponding to vertices of our
graph $\Ga$. Certainly, every trivalent graph $\Ga$ is geometrically
equivalent to the handlebody $\wGa$ or to the Riemann surface
$\Si_{\Ga}$ with fixed pants decomposition (4.1). Moreover, we have the
identification
\begin{equation}
 \pi_1(\wGa)=\pi_1(\Ga)
\end{equation}
and the epimorphism
\begin{equation}
r\colon \pi_1(\Si_{\Ga})\to \pi_1(\wGa)=\pi_1(\Ga).
\end{equation}

We can consider our graph $\Ga$ as a skeleton inside the handlebody
$\wGa$. Then for every edge $e$ inside the full tube with the boundary
$\we$ there exists a disc $\De$ in the handlebody with the boundary on
$\we$ meeting the edge $e$ transversely at one point such that
\begin{equation}
 \p \De=C_e,
\end{equation}
and we can move $\Ga$ to the boundary of the handlebody to get an
embedding
\begin{equation}
 j\colon\Ga\to \Si_{\Ga}
\end{equation}
defined up to isotopy.

Now fix a point $p_e\in C_e$ on each circle $C_e$ and a vertex $v_0\in
V(\Ga)$. We identify it with the point $j(v_0)\in \Si_{\Ga}$. Then we can
view (4.3) as being the epimorphism of the exact sequence
\begin{equation}
1\to \ker r\to \pi_1(\Si_{\Ga}, v_0) \xrightarrow{\,r\,}
\pi_1(\Ga,
 v_0)\to 1
\end{equation}
if we fix paths from $v_0$ to every point $p_e$.

The orientation of every path of length 1 $\vec e\in P_1(\Ga)$ and the
pumping normal vector field define an orientation of the tube $\we$ and
hence an orientation of the meridian $C_e$ of (4.1). Thus we get the cycle
$C_{\vec e}$ joined by the path with the point $v$. This gives a
collection of homotopy classes considered as elements of the fundamental
group $\pi_1 (\Si_\Ga, v_0)$
\begin{equation}
 \bigl\{C_{\vec e} \bigm| e\in E(\Ga) \bigr\}\subset \pi_1(\Si_{\Ga},v_0).
\end{equation}

More precisely,
\begin{equation}
\{ C_{\vec e} \} \subset \ker r.
\end{equation}

Now sending $\vec e$ to $ C_{\vec e} $ and a loop of length $d$ in $\Ga$
given by a sequence $\vec e_1,\dots, \vec e_d$ to
\begin{equation}
p(\vec e_1,\dots, \vec e_d)=C_{\vec e_d} \circ\cdots\circ C_{\vec
e_1}\in \pi_1(\Si_{\Ga}, v_0)
\end{equation}
gives a map of the collection of loops through $v_0$ to
$\pi_1(\Si_{\Ga}, v_0)$.
Hence we get the homomorphism
\begin{equation}
 q\colon \pi_1(\Ga, v_0)\to \ker r \subset \pi_1(\Si_{\Ga}, v_0).
\end{equation}
The composite $r \circ q$ defines the ``complex'' of groups
\begin{equation}
 \begin{matrix}
\cdots \xrightarrow{r \circ q} \pi_1(\Si_{\Ga}, v_0) \xrightarrow{r \circ
q} \pi_1(\Si_{\Ga}, v_0) \xrightarrow{r \circ q}\cdots \\[8pt]
(r\circ q)^2=\id
 \end{matrix}
\end{equation}
with cohomology
\begin{equation}
\ker(r \circ q) / \im (r \circ q)=\ker r / \im q=H(\Ga).
\end{equation}

%(COMPUTE FOR $g=2$)

\subsection{From gauge theory on $\Si_{\Ga}$ to gauge theory on $\Ga$}

Let $\sA_F (\Si_{\Ga})$ be the space of $\SU(2)$-connections on the
trivial vector bundle of rank 2 on $\Si_{\Ga}$, $\sG$ the corresponding
gauge group and
 \begin{equation}
\Rep(\pi_1(\Si_{\Ga}, v_0)=\sA_F(\Si_{\Ga}) /\sG_{v_0}
 \end{equation}
 the framed gauge-orbit space of flat connections.

The main construction relating gauge theories on a Riemann surface
$\Si_{\Ga}$ and on $\Ga$ is the map
\begin{gather}
m\colon\Rep(\pi_1(\Si_{\Ga}, v_0))\to
 \sA(\Ga) \\
m(\rho)(\vec e)=\rho (C_{\vec e}). \nonumber
\end{gather}
$\SU(2)$ acts on the target space of this map by the adjoint action and
on the source as the diagonal subgroup $\SU(2)_\De$ under the action of
the gauge group $\sG$. This map is equivariant with respect to these
actions, so we can factorize it as a map of quotients
\begin{equation}
m\colon \CLRep(\pi_1(\Si_{\Ga}))\to
 \CL\sA(\Ga)
\end{equation}
(see (3.16)). For any sequence $e_1,e_2,\dots,e_n$ of edges, we
construct the map
\begin{equation}
m_{(e_1, e_2,\dots , e_n)}=\prod_{i=1}^n i_{e_i} \circ m,
\end{equation}
where $i_{e_i}$ are involutions (2.4) and (3.2).

Different composites of this map with maps from gauge theory on $\Ga$
provide important fibrations of $\CLRep(\pi_1(\Si_{\Ga}))$. The first of
these is the composite
\begin{equation}
\conj \circ m=\pi_\Ga\colon \CLRep(\pi_1(\Si_{\Ga}))\to \De_\Ga
\subset \prod_{e\in E(\Ga)} [0, 1]_e,
\end{equation}
(see (3.28)) which is nothing other than the real polarization (see
\cite{6}, \cite{7} and below). This fibration is given by the moment map
of $(3g-3)$-torus action and is described perfectly.

The second fibration is given by the quotient projection
\begin{equation}
P \circ m\colon \Rep(\pi_1(\Si_{\Ga}))\to \Rep(\pi_1(\Ga))).
\end{equation}
From the direct description of maps it is easy to prove the following

\begin{prop} This map is nothing other than the composite
\begin{equation}
\Rep(\pi_1(\Si_{\Ga}, v_0)) \xrightarrow{\mathrm{res}}\Rep(\ker\quad
r, v_0) \xrightarrow{q^*}\Rep(\pi_1(\Ga)), v_0)
\end{equation}
where the first map is the restriction map to a subgroup and the second
is induced by the group homomorphism $q$ (4.10).
\end{prop}
Recall that $\CLRep(\pi_1(\Ga))$ is called the unitary Schottky
space of genus $g$ and denoted by the symbol $ uS_{g}$ \cite{5}.
This space can be presented as a homogeneous space
 \begin{equation}
 uS_{g}=\SU(2)^{g} / \Ad_{\mathrm{diag}} \SU(2)
 \end{equation}
where $\Ad_{\mathrm{diag}}\SU(2)$ is the diagonal adjoint action on the
direct product (see \cite{5}). The map $r^*$ (4.3) embeds the unitary
Schottky space
\begin{equation}
r^*\colon uS_{g} \hookrightarrow \CLRep(\pi_1(\Si_{\Ga}))
\end{equation}
but on the other hand the composite (4.18) gives the surjective
map
\begin{equation}
P_{cl} \circ m\colon \CLRep(\pi_1(\Si_{\Ga}))\to uS_{g}.
\end{equation}
Thus we have the map
\begin{equation}
P_{cl} \circ m \circ r^*\colon \CLRep(\pi_1(\Si_{\Ga}))\to
\CLRep(\pi_1(\Si_{\Ga}))
\end{equation}
inducing the cohomology homomorphism
\begin{equation}
(P_{cl} \circ m \circ r^*)^*\colon H^*(\CLRep(\pi_1(\Si_{\Ga})),\Z)\to
H^*(\CLRep(\pi_1(\Si_{\Ga})),\Z).
\end{equation}
The Florentino conjecture compares the target spaces of two projections
$P\circ m $ and $\conj\circ m=\pi_\Ga$, that is, the spaces $uS_{g}$ and
$\De_\Ga$.
\begin{con}[C. Florentino] There exists a one-to-one map
 \[
f\colon \De_\Ga\to uS_{g}
 \]
which is smooth inside $\De_\Ga$ and on the nonsingular part of $uS_{g}$
and continuous everywhere.
\end{con}

In the genus 2 case
\begin{equation}
 \De_\Theta=\CLRep(\pi_1(\wv_i))=uS_{2}
\end{equation}
(see \cite[Section~3]{6}). In what follows, we
\begin{enumerate}
\item give $uS_{g}$ the structure of a polytope (more precisely, the
direct product of $g-1$ copies of 3-dimensional tetrahedrons ) and
\item construct the natural polytope embedding $\De_\Ga$ to $uS_{g}$.
\end{enumerate}

 \subsection{From gauge theory on $\Si_{\Ga}$ to gauge theory on a trinion}

 Recall from (3.16) that $\CL \sA(\Ga)=\sA(\Ga) / \SU(2)_\De$. For every
trinion $\wv$ we have the map
 \begin{equation}
 m_v\colon \CLRep(\pi_1(\wv))\to
 \CL\sA(\Ga).
 \end{equation}
 \[
 m_v(\rho)(\vec e)=\rho (C_{\vec e}).
 \]
 for $e\in E(\Ga)_v$.

 Every vertex $v$ has three local edges $ ( e_1, e_2, e_3)=E(\Ga)_v$ which
define the 3-cube
 \begin{equation}
 [0,1]^{3}_v=\prod_{e\in E(\Ga)_v} [0,1]_{e}.
 \end{equation}

 \begin{prop}\label{pro!4.2}
 \begin{enumerate}
 \item The image of the composite
 \[
 m_v\circ \conj\colon \CLRep(\pi_1(\wv))\to [0,1]^{3g-3}=\prod_{e\in E(\Ga)}
 [0,1]_e
 \]
 is contained in $[0,1]^{3}_v=\prod_{e\in E(\Ga)_v}
 [0,1]_e$.
 \item The image $m (\CLRep(\pi_1(\wv)))$ is a section of
 the projection $\conj$ (3.28)
 \begin{equation}
 m_v\circ \conj\colon \CLRep(\pi_1(\wv))\to \prod_{v\in
 E(\Ga)_v} [0,1]_{e}.
 \end{equation}
 \item The image $\conj(\CLRep(\pi_1(\wv)) )$ is
 the tetrahedron $\De_v $ inscribed in the cube given by triangle
 inequalities
 \[
 \vert t_1-t_2 \vert \leq t_3 \leq t_1-t_2
 \]
 where $t_i$ is a coordinate on the interval $[0,1]_{e_i}$ for $e_i\in
 E(\Ga)_v$. Hence
 \item
 \[
 \CLRep(\pi_1(\wv))=\De_v=\De_\Theta=uS_2.
 \]
 \end{enumerate}
 \end{prop}
 (for the proof see \cite[Proposition~3.1]{6}).

 Now for every $v\in V(\Ga)$ and its corresponding trinion $\wv$ we have
 the restriction map
 \begin{equation}
 \res_v\colon \CLRep(\pi_1(\Si_{\Ga}))\to
 \CLRep(\pi_1(\wv)).
 \end{equation}
 This map is induced by the natural map
 \begin{equation}
 r_v\colon \pi_1(\wv)\to \pi_1(\Si_{\Ga}).
 \end{equation}
 The group $\pi_1(\wv)$ is the free group with 2 generators.
 Let
 \begin{equation}
 I_v \subset \pi_1(\Si_{\Ga})
 \end{equation}
 be the image of $r_v$. Then
 \begin{equation}
 \im \res_v=\CLRep (I_v) \subset \CLRep(\pi_1(\wv))=\De_v
 \end{equation}
 (see Proposition~\ref{pro!4.2}, (4)). For example, if $v\in v(L(\Ga))$
is a vertex on a loop $e\in L(\Ga)$then $I_v=\Z$ and $\im
\res_v=[0,1]_{e}$.

 \subsection{Half Riemann surface $\Si_\Ga$}

 We now give the hyperbolic graph (2.23) a one-to-one map
 \begin{equation}
 h\colon V_+\to V_-.
 \end{equation}
Then as in Figure~2, p.~\pageref{fig!1},
 \begin{figure}[htb]
 \centerline{\epsfig{file=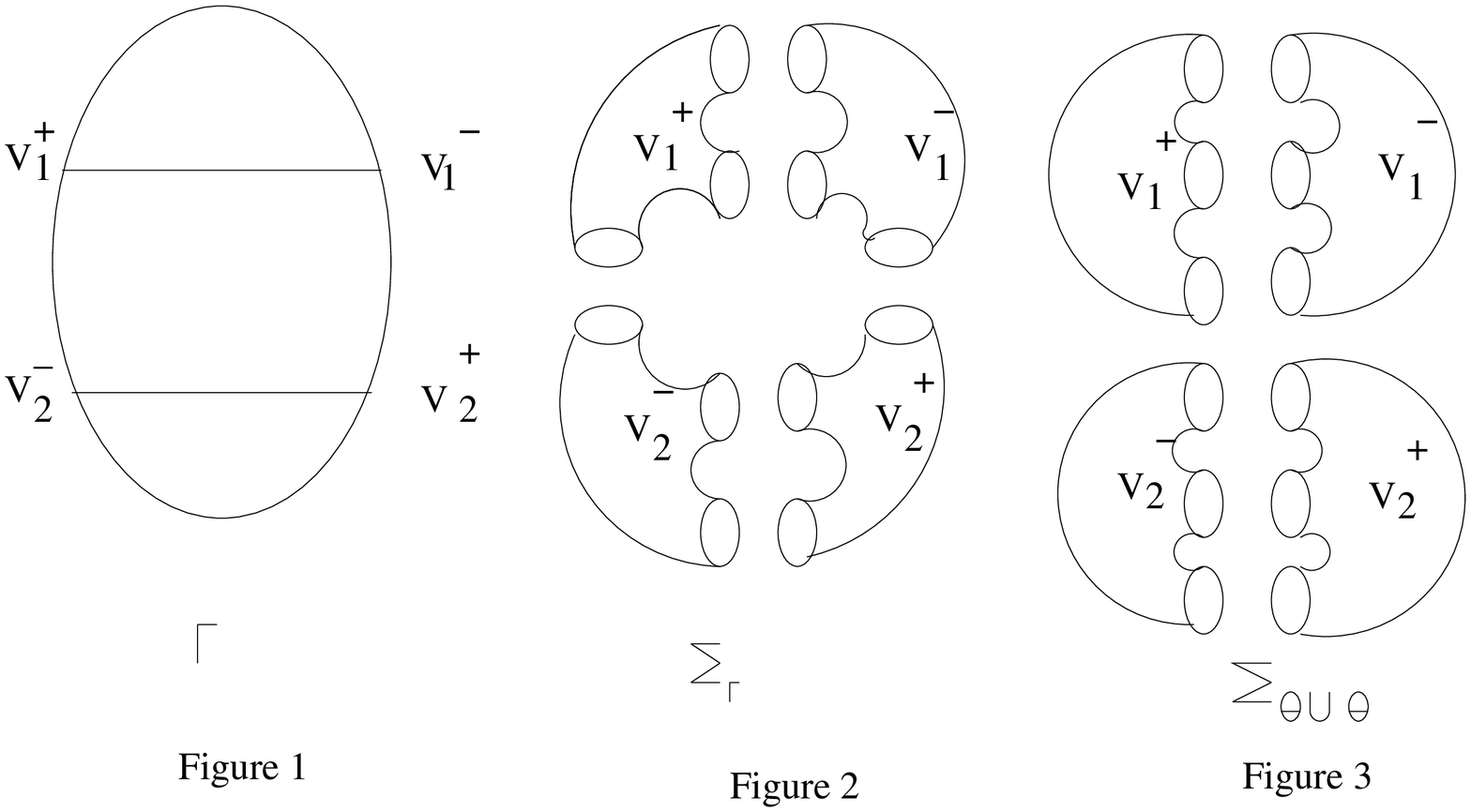,height=9cm,width=13cm}}
 \label{fig!1}
 \end{figure}
 we have:
 \begin{enumerate}
 \item
 \[
 V_+ \cap h(V_+)=\emptyset,\quad (V_+ \cup V_-)= V(\Ga),
 \quad
 \vert V_\pm \vert=g-1;
 \]
 \item for every pair of distinct vertices $v_i, v_j\in V_+$,
 \[
 E(\Ga)_{v_i} \cap E(\Ga)_{v_j}=\emptyset.
 \]
 \end{enumerate}

 \begin{dfn} The union
 \[
 \Si_\Ga^+=\bigcup_{v\in V_+} \wv
 \]
 is called a half of the Riemann surface $\Si_\Ga$ .
 \end{dfn}
 Thus a half of the Riemann surface $\Si_\Ga^+$ is the disjoint union of
trinions and it has a mirror $\Si_{\Ga}^{-}=h(\Si_\Ga^+)$ such that
 \begin{equation}
 \Si_\Ga^+ \cup \Si_\Ga^- =\Si_\Ga.
 \end{equation}
 A decomposition of this type is an analogue of a Heegard decomposition in
the 2-dimensional case and we can apply to this situation all the
constructions of \cite[Seminar~4]{1} (see ``Verlinde property''). The two
halves $\Si_\Ga^+$ and $\Si_\Ga^-$ are completely symmetric, and gluing is
given by the matrix $q_{\Ga}$ (see the end of section 2). Thus this matrix
plays the role of an element of the modular group in the 3-dimensional
case.

Now the fundamental group of $\Si_\Ga^+$ is given by
 \begin{equation}
 \pi_1(\Si_\Ga^+)=\prod_{v\in V_+} \pi_1(\wv),
 \end{equation}
 and hence
 \begin{equation}
 \CLRep(\pi_1(\Si_\Ga^+))=\prod_{v\in V_+} \CLRep(\pi_1(\wv)).
 \end{equation}
 Thus we have a natural homomorphism
 \begin{equation}
 r_+ (\Ga)\colon \pi_1(\Si_\Ga^+)\to \pi_1(\Si_\Ga)
 \end{equation}
 and the induced map
 \begin{equation}
 r_+^*(\Ga)\colon \CLRep (\pi_1(\Si_\Ga))\to \CLRep
(\pi_1(\Si_\Ga^+))=\prod_{v\in V_+} \CLRep(\pi_1(\wv)).
 \end{equation}

 The boundary of our half of the Riemann surface $\Si_\Ga^+$
 \begin{equation}
 \p \Si_\Ga^+=\{ C_e \} , \quad e\in E(\Ga)
 \end{equation}
 is the disjoint union of cycles (4.1) with the natural orientation coming
from the orientations of trinions.

 \subsection{Half Riemann surface $\Si_{\Theta^{g-1}}$}

 For genus 2 there is a unique connected trivalent graph without loops
$\Theta$ (denoted by this symbol for the shape). This graph gives the
disconnected graph of genus $g$
 \begin{equation}
 \Theta^{g-1}=\Theta \cup \cdots \cup \Theta
 \end{equation}
 which is disjoint union of $g-1$ copies of $\Theta$.

 In particular the fundamental group
 \begin{equation}
 \pi_1(\Theta^{g-1} )=\prod_{\text{$g-1$ copies}} \pi_1(\Theta)
 \end{equation}
 We can apply all the constructions of Section~3.1 to define the space of
connections $\sA(\Theta^{g-1})$, the gauge group $\sG(\Theta^{g-1}) $
action and the gauge orbit space of Section~3.3
 \begin{equation}
 \sA(\Theta^{g-1}) / \sG(\Theta^{g-1})=\CLRep (\pi_1(\Theta^{g-1}))
 \end{equation}
 and the map $\conj$ of Section~3.4.

 Moreover pumping $\Theta$ we get the compact Riemann surface $\Si_\Theta$
of genus 2 with the trinion decomposition of Figure~3. Thus the Riemann
surface $\Si_{\Theta^{g-1}}$ is the disjoint union
 \begin{equation}
 \Si_{\Theta^{g-1}}=\Si_\Theta \cup \cdots \cup \Si_\Theta
 \end{equation}
 of $g-1$ copies of a Riemann surface of genus 2 with the trinion
decomposition (4.1).

 A small difference between the case of connected and disjoint curves is
that we can't fix a basic point $v_0$. Thus we have to consider spaces of
representation classes. We have the map (4.15)
 \begin{equation}
 m\colon \CLRep(\pi_1(\Si_{\Theta^{g-1}}))\to
 \CL\sA(\Theta^{g-1}).
 \end{equation}
 and the composite (4.7)
 \begin{equation}
 \conj\circ m=\pi_{\Theta^{g-1}}\colon
 \CLRep(\pi_1(\Si_{\Theta^{g-1}}))\to \De_{\Theta^{g-1}}
 \subset \prod_{e\in E(\Theta^{g-1})} [0, 1]_e.
 \end{equation}

 Now consider the decomposition of $\Si_\Ga$ as an union of disjoint
trinion as a many piece jigsaw puzzle (as in Figure~2). We can transform
this collection to the union of disjoint trinions of a decomposition
of $\Si_{\Theta^{g-1}}$ (as in Figure~3). Thus we may consider the half of
the Riemann surface $\Si_\Ga $ as the half of the Riemann surface
$\Si_{\Theta^{g-1}}$ with the decomposition
 \[
 V_+ \cap h(V_+)=\emptyset,\quad V_+ \cup h( V_+)=
 V(\Theta^{g-1}), \quad
 \vert V_\pm \vert=g-1.
 \]
 Thus
the union $\Si_{\Ga}^{+}$ is a half of two curves
 \begin{equation}
 \Si_\Ga \supset \Si_{\Ga}^{+}
\subset \Si_{\Theta^{g-1}}
 \end{equation}
 and together with the homomorphism $r_+(\Ga)$ (4.38), we get another
 \begin{equation}
 r_+ (\Theta^{g-1})\colon \pi_1(\Si_\Ga^+)=\pi_1(\Si_{\Theta^{g-1}}^+)\to
 \pi_1(\Si_{\Theta^{g-1}})
 \end{equation}
 which componentwise is the standard map for Riemann surfaces of genus 2.

 \section{Gauge quotient map}

 \subsection{The section of $\conj$}

 We have identifications
 \[
 E(\Ga)=E(\Theta^{g-1}), \quad V(\Ga)=V(\Theta^{g-1}), \quad
P_1(\Ga)=P_1(\Theta^{g-1})
 \]
 and hence
 \begin{equation}
 \sA(\Ga)=\sA(\Theta^{g-1})
\quad\hbox{and}\quad \CL\sA(\Ga)=\CL\sA(\Theta^{g-1}).
 \end{equation}
 But the gauge groups $\sG(\Ga)$ and $\sG(\Theta^{g-1})$ actions are
entirely different. In particular, the gauge quotient map for the gauge
group $\sG(\Theta^{g-1})$ can be decomposed as follows
 \begin{equation}
 \CL\sA(\Ga) \xrightarrow{p_\Theta^{g-1}} \prod_{(g-1) \text{ copies }}
 \CL\sA(\Theta)
 \xrightarrow{(g-1)\times p_\Theta} \prod_{(g-1) \text{
 copies }}\CLRep(\pi_1(\Theta)).
 \end{equation}

 All components of the latter map are the same
 \begin{equation}
 p_\Theta\colon \CL\sA(\Theta)\to \CLRep(\pi_1(\Theta)).
 \end{equation}
 We study this map later. The maps $\conj$ are absolutely the same, with
the same target space $\prod_{e\in E(\Ga)}[0, 1]_e$. Moreover, the product
of the maps $m_v$ (4.26) gives the map
 \begin{equation}
 m_+=\prod_{v\in V_+} m_v\colon \prod_{v\in V_+}
 \CLRep(\pi_1(\wv))\to \CL \sA(\Ga).
 \end{equation}
 Let
 \begin{equation}
 W=m_+ \Bigl(\prod_{v\in V_+} \CLRep(\pi_1(\wv))\Bigr) \subset \CL
\sA(\Ga)
 \end{equation}
 be the image of (5.4). A direct interpretation of the maps in the
composite
 \begin{equation}
 \CLRep (\pi_1(\Si_\Ga)) \xrightarrow{r^*_+(\Ga)} \prod_{v\in V_+}
 \CLRep(\pi_1(\wv)) \xrightarrow{m_+} \CL \sA(\Ga),
 \end{equation}
 gives

 \begin{prop}
 \begin{enumerate}
 \item The image
 \[
 r^*_+ (\Ga)\circ m_+ (\CLRep (\pi_1(\Si_\Ga))) \subset W.
 \]
 \item The map (4.17)
 \[
 \conj\circ m\colon \CLRep (\pi_1(\Si_\Ga))\to \De_\Ga \subset
\prod_{e\in E(\Ga)} [0, 1]_e
 \]
 is the composite of $r^*_+(\Ga)$ and
 \begin{equation}
 \conj\circ m_+\colon \prod_{v\in V_+} \CLRep(\pi_1(\wv))\to
 \prod_{v\in V_+} \De_v \subset \prod_{e\in E(\Ga)} [0, 1]_e.
 \end{equation}
Thus
 \item
 \begin{equation}
 \De_\Ga \subset \prod_{v\in V_+} \De_v \subset \prod_{e\in E(\Ga)} [0,
 1]_e.
 \end{equation}
 \end{enumerate}
 \end{prop}

 \begin{prop}
 The image $W$ is a section of the projection $\conj$. That is, the
restriction of $\conj$
 \[
 \conj \rest{W}\colon W\to \prod_{v\in V_+} \De_v \subset \prod_{e\in
 E(\Ga)} [0, 1]_e
 \]
 is one-to-one.
 \end{prop}
 To prove this statement, we may use the $\conj$--map for $\Theta^{g-1}$.
Recall that trinions of components of the direct product $W$ don't have
boundaries in  common. Thus coordinates of the cube are divided in triples
corresponding vertices $v\in V_+(\Ga)=V_+(\Theta^{g-1})$. Then 
Proposition~4.2, (3) implies the required statement componentwise.

 \begin{cor} The identification $W=\prod_{v\in V_+} \De_v $
 determines the section
 \begin{equation}
 s\colon \prod_{v\in V_+} \De_v\to \CL \sA(\Ga)
 \end{equation}
 and the restriction of this section to $\De_\Ga \subset \prod_{v\in V_+}
 \De_v$
 determines the embedding
 \begin{equation}
 s \rest{\De_\Ga}\colon \De_\Ga\to \CL\sA(\Ga).
 \end{equation}
 \end{cor}

C. Florentino gave a beautiful description of the subpolytope $\De_{\Ga}
\subset \De_{\Theta}^{g-1}$ in terms of symmetries of a graph $\Ga$.

Now the restriction $s \rest{\De_\Ga}(\De_\Ga)$ of the quotient map of
$\sG(\Ga)$-action sends $\De_\Ga $ to $uS_g$. Hence it is just what we need
to construct the natural map $f\colon \De_\Ga\to uS_g$ for the Florentino
conjecture. But the key is the restriction of the quotient map of
$\sG(\Ga)$-action to $W$, that is, comparing the gauge theory on $\Ga$
and the gauge theory on $\Theta^{g-1}$:
 \begin{equation}
 s\circ P_{cl}\colon W\to \CLRep (\pi_1(\Ga)),
 \end{equation}
 where $P_{cl}$ is the projection of (3.22).

 \subsection{$\Ga=\Theta$}

 In this case the handlebody $\wTh$ is a full pretzel. We get the Riemann
surface $\Si_{\Theta}$ of genus $8$ marked by tubes $\we_1,\we_2,\we_6$
for every edge of the triple $E(\Ga)$ and trinions $\wv_+ $ and $
\wv_- $ of the pair $ V(\Theta)$. The isotopy classes of
 meridian circles of tubes define 3 disjoint, non\-contractible,
 pairwise nonisotopic classes $c_i$ around $e_i$. The complement
 is the union
 \begin{equation}
 \Si_{\Theta} \setminus \{c_1,c_{2},c_{3}\}= \wv_{+}
 \cup \wv_{-}
 \end{equation}
 of 2 trinions. We fix a point $ p\in \wv_{+}$ and the natural orientation
so that
 \begin{equation}
 c_1 \cdot c_2 \cdot c_2=1
 \end{equation}
 as classes of pointed loops.

 Now
\begin{equation}
\sG(\Theta)=\SU(2)^{V(\Ga)}=\SU(2)_{+} \times \SU(2)_{-}
 \end{equation}
 with the action
 \begin{equation}
 g_{+} \cdot * \cdot g_{-}^{-1}
 \end{equation}
if $v_{+}=v_{s}$. We finally get the map
 \begin{equation}
 P_{cl}\colon \CLRep(\pi_1(\Si_{\Theta}))\to
 \CLRep(\pi_1(\Theta)).
 \end{equation}

 \begin{prop} The image $i_{e_{3}}(W)$ is a section of
 the projection $P_{cl}$ (3.22).
 \end{prop}
 The  boundary components $\p \wv=c_1, c_2, c_3$ as elements of the
fundamental group of $\wv$ constrain the equation (5.13). Consider the map
 \begin{equation}
 m_{e_3}\colon \CLRep (\pi_1(\wv))\to \CL \sA(\Theta)
 \end{equation}
 Then the component $g_v$ of an element of the gauge group acts on these
cycles as $g(\rho(c_i))=\rho(c_i)\circ g$ for $i =1,2$ and
$g(\rho(c_3))=g^{-1} \rho(c_3)$. To preserve the equation (5.13) we have
 \begin{equation}
 \rho(c_1)\circ g\circ \rho(c_2)\circ g\circ g^{-1} \rho(c_3)=1.
 \end{equation}
 From this we have $g=1$ and we are done.

 \begin{cor} This identification
 \begin{equation}
 W=\CLRep (\pi_1(\Theta))=uS_2
 \end{equation}
 defines
 \begin{enumerate}
 \item the polytope structure on $uS_2$ and
 \item the polytope equivalence
 \begin{equation}
 f\colon \De_\Theta\to uS_2.
 \end{equation}
 \end{enumerate}
 \end{cor}
Thus the Florentino conjecture is true in the case $g=2$. We used this
identification for the construction of non-Abelian theta functions in the
first nontrivial case $g=2$ considered in \cite{9}.

 \subsection{General case}

We show how to prove the statement of Proposition 5.3 for a multi-theta
graph of any genus such as Figure~3, p~\pageref{fig!1}. The principle of
the proof is just as before. Namely in (5.18) the element $g$ after
$\rho(c_1)$ doesn't have to be the same as after $\rho(c_{2})$. Then in
the same vein $g$ has to be 1.

 To apply this hint, remark that the set of vertices $V(\Ga)$ decompose in
pairs $(v_{i}^{+}, h(v_{i}^{+})= v_{i}^{-})$. Let us start with the pair
$v_1^{\pm}$. Using the adjoint action we may suppose that
$g_{v_1}^{-}=1$. Assume
\begin{equation}
\p \wv_1^{+} \cap \p \wv_1^{-}=c_{2}
\cup c_{3}
\end{equation}
and
 \[
c_1 =\p \wv_1^{+} \cap \p \wv_{2}^{-}
 \]
 We can choose the map $m_{\{e\}}$ of (4.16) so that $e_1$ and $e_{2}$
have opposite orientations. Then we have equations
 \[
c_1\circ c_{2}\circ c_{3}=1
 \]
and
\begin{equation}
\rho(c_1)\circ g_{2}^{-}\circ \rho(c_{2})\circ g_1^{+}\circ
(g_1^{+})^{-1} \circ \rho(c_{3})=1.
\end{equation}
Thus $g_{2}^{-}=1$. Applying the same construction to $\wv_1^{-}$ we get
$g_{2}^{+}=1$ also. This gives the proof for genus 3 case.

Iterating these arguments gives the proof for any hyperbolic graph of any
genus.

\section{Acknowledgments} I would like to express my gratitude to my
collaborators C.~Florentino, J.~Mourao, J.P.~Nunes and the Instituto
Superior Tecnico of Lisbon for support and hospitality. I also thank Miles
Reid for help with English.


\begin{thebibliography}{99}

\bibitem{1} Michael Atiyah, ``The path integral formulation'', Oxford
seminar on Jones--Witten theory, Michaelmas Term 1988. Seminar 6, p.~106.

\bibitem{FQ} Daniel S. Freed and Frank Quinn, ``Chern--Simons theory with
finite gauge group'', hep-th/9111004 44~pp.

\bibitem{2} A. N. Tyurin ``Special Langrangian geometry and slightly
deformed algebraic geometry (Splag and Sdag).'' Warwick Preprint: 22/1998,
section 1. math.AG/9806006, 45~pp.

\bibitem{3} N. Hitchin and J. Sawon, ``Curvature and characteristic
numbers of hyper\-K\"ahler manifolds'', math.DG/9908114, 17 pp.

 \bibitem{4} S. Donaldson ``Gluing techniques in the cohomology of
moduli spaces'', in Topological methods in modern mathematics (Stony
Brook, NY, 1991), Publish or Perish, Houston, TX, 1993, pp. 137--170.

 \bibitem{5} A. Tyurin, ``Three mathematical facets of SU(2)-spin
networks'',\newline math.DG/0011035, 20~pp.

 \bibitem{6} L. C. Jeffrey and J. Weitsman, Bohr--Sommerfeld orbits in
the moduli space of flat connections and the Verlinde dimension formula,
Commun. Math. Phys. {\bf150} (1992) 593--630

 \bibitem{7} A. Tyurin, Quantization and theta functions,
math.AG/9904046; On Bohr-Sommerfeld bases, math.AG/9909084; Izv. Ross.
Akad. Nauk Ser. Mat. {\bf64} (2000), 163--196; English transl. in Russian
Acad. Sci. Izv. Math.{\bf64}:5 (2000)

 \bibitem{8} T. R. Ramadas, L. M. Singer and J. Weitsman, Some comments
on Chern--Simons gauge theory, Commun. Math. Phys. {\bf126} (1989) 409--420

\bibitem{9} C. Florentino, J. Mourao, J. Nunes and A. Tyurin 
``Analytical aspects of the theory non-Abelian theta functions.'' Preprint
IST 2000.

 \end{thebibliography}
 \end{document}